\documentclass[a4paper, 12pt]{article}

\usepackage{amsfonts, amsmath, amsthm}
\usepackage{amssymb}
\usepackage{latexsym}
\usepackage[dvips]{graphicx}
\usepackage{color}

\theoremstyle{plain}
\newtheorem{thm}{Theorem}[section]
\newtheorem{prop}[thm]{Proposition}
\newtheorem{lem}[thm]{Lemma}
\newtheorem{cor}[thm]{Corollary}
\theoremstyle{definition}
\newtheorem{dfn}[thm]{Definition}

\theoremstyle{definition}
\newtheorem{rmk}[thm]{Remark}
\numberwithin{equation}{section}

\newcommand{\R}{\mathbb{R}}
\newcommand{\C}{\mathbb{C}}

\newcommand{\cc}[1]{\overline{#1}}
\newcommand{\op}[1]{\mathcal{#1}}

\newcommand{\pa}{\partial}
\newcommand{\eps}{\varepsilon}
\newcommand{\jb}[1]{\langle #1 \rangle}

\newcommand{\realpart}{{\rm Re\,}}
\newcommand{\imagpart}{{\rm Im\,}}

\newcommand{\dis}{\displaystyle}

\newcommand{\sh}[1]{{#1}^{\sharp}}
\newcommand{\tm}{\tilde{m}}

\begin{document}
\title{
Recent advances  on Schr\"odinger equations with 
dissipative nonlinearities}

\author{
             Chunhua Li  \thanks{
              Department of Mathematics, College of Science,
              Yanbian University.
              977 Gongyuan Road, Yanji, Jilin 133002, China.
              (E-mail: {\tt sxlch@ybu.edu.cn})
             }
          \and
          Yoshinori Nishii\thanks{
             Department of Mathematics, 
              Tokyo University of Science. 
              1-3 Kagurazaka, Shinjuku-ku, Tokyo 162-8601, Japan.
              (E-mail: {\tt yoshinori.nishii@rs.tus.ac.jp}) 
             }
          \and
          Yuji Sagawa \thanks{
              Mathematical Science Education Center,
              Kumamoto University. 
              2-40-1 Kurokami, Chuo-ku, Kumamoto 860-8555, Japan. 
             }
          \and
          Hideaki Sunagawa \thanks{
              Department of Mathematics, Graduate School of Science, 
              Osaka Metropolitan  University.
              3-3-138 Sugimoto, Sumiyoshi-ku, Osaka 558-8585, Japan. 
              (E-mail: {\tt sunagawa@omu.ac.jp})
             }
}

\date{\today }
\maketitle

\vspace{-4mm}
\begin{center}
Dedicated to Professor Tohru Ozawa on the occasion of
his sixtieth birthday\\
\end{center}
\vspace{2mm}

\noindent{\bf Abstract:}\ 
We give a survey on recent developments  on  nonlinear 
Schr\"odinger equations with dissipative structure   
based on the authors' recent works.




\section{Introduction} \label{sec_intro}

This paper is intended to be a survey on recent advances on nonlinear 
Schr\"odinger equations with dissipative structure based on the authors' 
recent works 
\cite{LNSS1}, \cite{LNSS2}, \cite{LNSS3}, \cite{LNSS4}. 
We refer the readers to these papers and the 
references cited therein for the detailed expositions. 

This paper is organized as follows. 
In Section \ref{sec_powerNLS}, we summarize typical previous 
results for nonlinear Schr\"odinger equations with the power-type 
nonlinearities. Section \ref{sec_DNLS} is 
devoted to the case where the nonlinear term depends also on the 
derivative of the unknown function. Special attentions are paid to 
the weakly dissipative nonlinearities which never appear 
in the power-type nonlinearity situation. 
In Section \ref{sec_two-compNLS}, we focus our attentions on a 
two-component Schr\"odinger system which tells us that the system case 
is much more delicate than the single case. 
Finally, we enumerate the results obtained in \cite{LiS1} concerning 
general nonlinear Schr\"odinger systems of derivative type in the Appendix. 
Throughout this paper, we denotes by $\op{L}$ the standard free 
Schr\"odinger operator $i\pa_t+\frac{1}{2}\pa_x^2$ for $(t,x)\in \R\times \R$ 
with $i=\sqrt{-1}$. 
The free evolution group $e^{\frac{it}{2}\pa_x^2}$ is written as $\op{U}(t)$. 
The function space $H^k$ stands for the $L^2$-based 
Sobolev space of order $k$ equipped with the norm 
$\|\phi\|_{H^k}=\sum_{0 \le j \le k}\|\pa_x^j \phi\|_{L^2}$,
and the weighted Sobolev space $H^{k,m}$ is defined by 
$\{\phi \in L^2\,|\, \jb{\, \cdot \,}^{m} \phi \in H^k \}$ with 
$\jb{x}=\sqrt{1+x^2}$. 
Several non-negative constants will be denoted by the same letter $C$, 
unless otherwise specified.

\section{Brief review on the basic facts}\label{sec_powerNLS} 
\label{sec_review}

First of all, let us recall some of well-known results 
on large-time behavior of small data solutions to the cubic power-type 
nonlinear Schr\"odinger equation in the form  
\begin{align}
 \op{L} u = \lambda |u|^{2} u, \qquad t>0,\ x\in \R,
 \label{nls_1}
\end{align}
where $\lambda$ is a constant. 
What is interesting in \eqref{nls_1} is 
that the large-time behavior of the solution is actually affected by the 
coefficient $\lambda$ even if the initial data is sufficiently small, 
smooth and decaying fast as $|x|\to \infty$. 
If $\lambda\in \R$, it is shown by Ozawa \cite{Oz} and Hayashi--Naumkin  
\cite{HN1} that the solution to \eqref{nls_1} with small data behaves 
like 
\[
 u(t,x)=\frac{1}{\sqrt{it}} \alpha(x/t) 
 e^{i\{\frac{x^2}{2t}  - \lambda |\alpha(x/t)|^2 \log t \}}
 +o(t^{-1/2})
\quad \mbox{as} \ \ t\to +\infty
\]
with a suitable $\C$-valued function $\alpha(y)$. 
An important consequence of this asymptotic expression is that the solution 
decays like $O(t^{-1/2})$ in $L^{\infty}(\R_{x})$, 
while it does not behave like free solutions unless $\lambda=0$. 
In other words, the additional logarithmic factor in the phase 
reflects the long-range character of the cubic nonlinear Schr\"odinger 
equations in one space dimension. 
If $\lambda\in \C\backslash \R$ in \eqref{nls_1}, 
another kind of long-range effect can be observed. 
For instance, according to \cite{Shim} 
(see also \cite{KitaShim}, \cite{JJL}, \cite{HLN}, etc.), 
the small data solution $u(t,x)$ to \eqref{nls_1} decays like 
$O(t^{-1/2}(\log t)^{-1/2})$ in $L^{\infty}(\R_x)$ as $t\to +\infty$ 
if $\imagpart\lambda<0$. 
This gain of additional logarithmic time decay should be interpreted as 
another kind of long-range effect 
(see also \cite{CH}, \cite{CHN}, \cite{HLN}, \cite{HLN2}, \cite{HNS}, 
\cite{Hoshino}, 
\cite{Hoshino2}, \cite{JJL}, \cite{KLS}, \cite{Kim}, \cite{KitaLi}, 
\cite{KitaNak}, \cite{KitaShim}, \cite{LNSS1}, \cite{LNSS2}, 
\cite{LiS1}, \cite{OgSat}, \cite{Sat1}, \cite{Sat2}, and so on). 
Time decay in $L^2$-norm is also investigated by several authors. 
Among others, it is pointed out by Kita-Sato \cite{KitaSato} that 
the optimal $L^2$-decay rate is $O((\log t)^{-1/2})$ in the case of 
\eqref{nls_1} with $\imagpart\lambda<0$. 
We are interested in extending $L^2$-decay results of this kind to 
derivative nonlinearity case or system case.

\section{Nonlinear Schr\"odinger equations of derivative type: 
Weak dissipativity}
\label{sec_DNLS}

In this section, we focus on the initial value problem in the form
\begin{eqnarray}
\op{L} u=N(u,\pa_x u),
\qquad  t>0,\ x\in \R
\label{eq_d}
\end{eqnarray}
with
\begin{eqnarray}
u(0,x)=\varphi(x), 
\qquad  x\in \R,
\label{data_sing}
\end{eqnarray}
where $\varphi$ is a prescribed $\C$-valued function on $\R$. 
The nonlinear term $N(u,\pa_x u)$ is a cubic homogeneous polynomial in  
$(u,\cc{u}, \pa_x u, \cc{\pa_x u})$ with complex coefficients. 
If $\varphi$ is $O(\eps)$ in $H^3\cap H^{2,1}$ with $0<\eps \ll 1$, 
what we can expect 
for general cubic nonlinear Schr\"odinger equations in $\R$ 
is the lower estimate for the lifespan $T_{\eps}$ 
in the form $T_{\eps}\ge \exp(c/\eps^2)$ with some $c>0$ not depending on 
$\eps$, 
and this is best possible in general (see \cite{Kita} for an example of 
small data blow-up). More precise information on the lifespan  
is available under the restriction 
\begin{align}
 N(e^{i\theta},0) =e^{i \theta} N(1,0), \qquad \theta \in \R
\label{weak_gi}
\end{align}
and the initial condition 
 \begin{eqnarray}
 u(0,x) =\eps \psi(x),\qquad x\in \R,
\label{data_eps_sing}
\end{eqnarray}
instead of \eqref{data_sing}, 
where $\psi\in H^{3}\cap H^{2,1}$ is independent of $\eps$. 
In fact we have the following. 
\begin{thm}[\cite{SagSu}, \cite{Su2}, \cite{Su3}] \label{thm_lifespan}
Assume that $\psi \in H^3 \cap H^{2,1}$. 
Suppose that the nonlinear term $N$ satisfies \eqref{weak_gi}. 
Let $T_{\eps}$ be the supremum of $T>0$ such that the 
initial value problem \eqref{eq_d}--\eqref{data_eps_sing} admits a 
unique solution in $C([0,T); H^3 \cap H^{2,1})$. 
Then it holds that 
\begin{align}
 \liminf_{\eps \to +0} \eps^2 \log T_{\eps} 
\ge 
\frac{1}{\dis{2\sup_{\xi \in \R}(|\hat{\psi}(\xi)|^2 \imagpart \nu(\xi))}}
\label{lifespan}
\end{align}
with the convention $1/0=+\infty$,
where the function $\nu:\R\to \C$ is defined by 
\begin{align}\label{def_nu}
 \nu(\xi)=\frac{1}{2\pi i} \oint_{|z|=1} N(z,i\xi z) \frac{dz}{z^2}, 
\quad \xi \in \R,
\end{align}
and $\hat{\psi}$ denotes the Fourier transform of $\psi$, i.e.,
\[
 \hat{\psi}(\xi)=\op{F}\psi(\xi)=
 \frac{1}{\sqrt{2\pi}} \int_{\R} e^{-iy\xi} \psi(y)\, dy, 
\quad \xi \in \R.
\]
\end{thm}
Note that \eqref{weak_gi} excludes just the worst terms $u^3$, $|u|^2\cc{u}$, 
$\cc{u}^3$. It is known that these three terms are quite difficult to 
handle in the present setting, and we do not pursue this case here 
(cf. \cite{MP}).

In view of the right-hand side in \eqref{lifespan}, it may be natural 
to expect that the sign of $\imagpart \nu (\xi)$ has something to do with 
global behavior of small data solutions to \eqref{eq_d}. 
In fact, it has been pointed out in \cite{SagSu} that typical results on 
small data global existence and large-time asymptotic behavior for 
\eqref{eq_d} under \eqref{weak_gi}  can be summarized in terms of 
$\imagpart \nu(\xi)$ as follows:
\begin{itemize}
\item
Small data global existence holds in $C([0,\infty);H^3\cap H^{2,1})$ 
under the condition
\begin{align}
\imagpart \nu(\xi) \le 0,\quad \xi \in \R.
\tag{${\bf A}$}
\end{align}
(See also Theorem~\ref{thm_sdge} in Appendix.)
\item
The global solution has (at most) logarithmic phase correction if 
\begin{align}
\imagpart \nu(\xi) = 0,\quad \xi \in \R.
\tag{${\bf A}_0$}
\end{align}
Also it is not difficult to see that  there is no $L^2$-decay under 
(${\bf A}_0$) for generic initial data of small amplitude.
\item
$L^2$-decay of the global solution occurs under the condition 
\begin{align}
\sup_{\xi \in \R} \imagpart \nu(\xi)<0.
\tag{${\bf A}_+$}
\end{align}
(See also Theorem~\ref{thm_decay1} in Appendix.)
\end{itemize}

Note that $\nu(\xi)=\lambda$ if $N=\lambda |u|^2u$. So these results 
cover the results in the power-type nonlinearity case mentioned in 
Section~\ref{sec_review}. However, as pointed out in \cite{LNSS3}, 
an interesting case is not covered by these classifications, 
that is the case where  $({\bf A})$ is satisfied but $({\bf A}_0)$ and 
$({\bf A}_+)$ are violated. For example, if $N=-i|\partial_x u|^2u$, 
we can easily check that $\imagpart \nu(\xi)=-\xi^2\le 0$, 
while the inequality is not strict because of vanishing at $\xi=0$. 
This is what we are interested in. 

To going further, let us remember the fact that, 
if $({\bf A})$ is satisfied but $({\bf A}_0)$ and $({\bf A}_+)$ are violated, 
then there exist $c_0>0$ and $\xi_0\in \R$ such that 
$\imagpart \nu(\xi)=-c_0(\xi-\xi_0)^2$. 
The converse is also true. 
This fact naturally leads us to the following definition of 
{\em the weak dissipativity}.

\begin{dfn} \label{dfn_wd}
We say that a cubic nonlinear term $N$  
is {\em{weakly dissipative}} if  the following two conditions (i) and 
(ii) are satisfied: 
\begin{itemize}
\item[(i)] $N(e^{i\theta},0) =e^{i \theta} N(1,0)$ for $\theta \in \R$.
\item[(ii)]  There exist $c_0>0$ and $\xi_0\in \R$ such that 
$\imagpart \nu(\xi)=-c_0(\xi-\xi_0)^2$.
\end{itemize}
\end{dfn}

The following two results reveal the $L^2$-decay property in the 
weakly dissipative case. 
\begin{thm}[\cite{LNSS4}]  \label{thm_upper}
Suppose that  $N$ is weakly dissipative 
and that $\eps =\|\varphi\|_{H^3\cap H^{2,1}}$ is sufficiently small. 
Then there exists a positive constant $C$, not depending on $\eps$, 
such that the global solution $u$ to \eqref{eq_d}--\eqref{data_sing} satisfies
\[
\|u(t)\|_{L_x^2}\le \frac{C\eps}{(1+ \eps^2\log (t+1))^{1/4}}
\] 
for $t\ge 0$.
\end{thm}

\begin{thm}[\cite{LNSS4}] \label{thm_lower}
Suppose that $N$ is weakly dissipative 
and  that the Fourier transform of $\psi$ does not 
vanish at the point $\xi_0$ coming from {\rm (ii)} in 
{\upshape{Definition~\ref{dfn_wd}}}.
Then we can choose  $\eps_0>0$ such that 
the global solution $u$ to \eqref{eq_d}--\eqref{data_eps_sing} 
satisfies
\[
 \liminf_{t\to +\infty}((\log t)^{1/4}\|u(t)\|_{L_x^{2}})>0
\]
for  $\eps \in (0,\eps_0]$.
\end{thm}

\begin{rmk}  According to \cite{KitaSato}, 
the optimal $L^2$-decay rate is $O((\log t)^{-1/2})$ in the case where
$N=\lambda |u|^2u$ with $\imagpart\lambda<0$. 
This should be contrasted with 
Theorems \ref{thm_upper} and \ref{thm_lower}, 
because these tell us that the optimal $L^2$-decay rate in the 
weakly dissipative case is $O((\log t)^{-1/4})$.

\end{rmk}

Now, let us explain heuristically why $L^2$-decay rate should be 
$O((\log t)^{-1/4})$  if $\hat{\psi}(\xi_0)\ne 0$. 
For this purpose, let us first remember the fact that the solution $u^0$ 
to the free Schr\"odinger equation (i.e., the case of $N=0$) 
behaves like 
\[
 \pa_x^k u^0(t,x)
 \sim 
 \left(\frac{ix}{t}\right)^k \frac{e^{-i\pi/4}}{\sqrt{t}} 
 \hat{\varphi}\left( \frac{x}{t}\right) e^{i\frac{x^2}{2t}}
 +\cdots
\]
as $t\to+\infty$ for $k=0,1,2,\ldots$. 
Viewing it as a rough approximation of the solution 
$u$ for \eqref{eq_d}, we may expect that $\pa_x^k u(t,x)$ could be 
better approximated by 
\[
 \left(\frac{ix}{t}\right)^k \frac{1}{\sqrt{t}} 
 A\left(\log t,  \frac{x}{t}\right) e^{i\frac{x^2}{2t}}
\]
with a suitable function $A(\tau,\xi)$, 
where $\tau=\log t$, $\xi=x/t$ and $t\gg 1$. 
Note that 
\[A(0,\xi)=e^{-i\pi/4}\, \hat{\varphi}(\xi)\] 
and that the extra variable $\tau=\log t$ is responsible 
for possible long-range nonlinear effect. 
Substituting the above expression into \eqref{eq_d}
and keeping only the leading terms, 
we can see (at least formally) that $A(\tau,\xi)$ should satisfy 
the ordinary differential equation 
\begin{align*}
 i\pa_{\tau} A = \nu(\xi)|A|^2A+\cdots 
\end{align*}
under \eqref{weak_gi}. If $N$ is weakly dissipative, we see that 
\[
  \pa_{\tau} |A|^2 = -2c_0(\xi-\xi_0)^2|A|^4+\cdots.
\]
Then it follows that 
\[
 |A(\tau,\xi)|^2
=
\frac{|\hat{\varphi}(\xi)|^2}{1+2c_0(\xi-\xi_0)^2 |\hat{\varphi}(\xi)|^2\tau} 
+\cdots,
\]
whence 
\begin{align*}
\|u(t)\|_{L_x^2}
\sim 
\|A(\log t)\|_{L^2_{\xi}}
\sim 
\left(\int_{\R} 
\frac{|\hat{\varphi}(\xi)|^2}{1+2c_0(\xi-\xi_0)^2|\hat{\varphi}(\xi)|^2
\log t}\, d\xi\right)^{1/2}\quad 
(t\to+\infty).
\end{align*}
By considering the behavior as $t\to +\infty$ of this integral carefully, 
we see that $L^2$-decay rate in the weakly dissipative case should be 
just $O((\log t)^{-1/4})$ if $\hat{\varphi}(\xi_0)\ne 0$.
Indeed, we have the following lemma. 
\begin{lem}[\cite{LNSS4}] \label{lem_int1}
Let $\theta \in L^{\infty}$ and $\xi_0\in \R$. We set 
\begin{align*}
S(\tau)
= \int_{\R} \frac{|\theta(\xi)|^2}{1+(\xi-\xi_0)^2|\theta(\xi)|^2\tau}\, d\xi
\end{align*}
for $\tau\ge 1.
$\begin{itemize}
\item[\upshape{(1)}]We have 
\[
 S(\tau) \le 4\|\theta\|_{L^{\infty}}\tau^{-1/2}, \qquad \tau \geq 1.
\]
\item[\upshape{(2)}]
Assume that there exists an open interval $I$ with $I\ni \xi_0$ such that 
$\inf_{\xi \in I}|\theta(\xi)|>0$.
Then we can choose a positive constant $C_*$, which is independent of 
$\tau \ge 1$ {\upshape{(}}but may depend on $\theta$ and 
$\xi_0${\upshape{)}}, 
such that 
\[
 S(\tau) \ge C_*\tau^{-1/2}, \qquad \tau \geq 1.
\]
\end{itemize}
\end{lem}

Our strategy of the proof of Theorems~\ref{thm_upper} and \ref{thm_lower} 
is to justify the above heuristic argument, which has been carried out in 
\cite{LNSS3} and \cite{LNSS4}. The key is to concentrate on the function 
\[
\alpha(t,\xi)=\op{F}[\op{U}(-t)u(t,\cdot)](\xi),
\] 
which is expected to play the role of $A(\log t,\xi)$ in the above argument. 
For the details, see \cite{LNSS3} and \cite{LNSS4}.

\section{A two-component system of nonlinear Schr\"odinger equations}
\label{sec_two-compNLS}

In this section, we turn our attentions to the system case.
Our goal is to reveal that the system case is much more 
delicate than the single case by considering the specific two-component system 
\begin{equation}
\left\{\begin{array}{ll}
\begin{array}{l}
 \op{L} u_1=- i|u_2|^2 u_1,\\
 \op{L} u_2=- i|u_1|^2 u_2,
 \end{array}
 & (t,x )\in (0,\infty)\times \R
 \end{array}\right.
\label{eq}
\end{equation}
under the initial condition 
\begin{equation}
u_j(0,x)=\varphi^0_j(x),
\qquad  x \in \R,\ j=1,2.
\label{data}
\end{equation}
When $\varphi_1^0=\varphi_2^0$, the system \eqref{eq}--\eqref{data} 
is reduced  to the single equation (\ref{nls_1}) with $\lambda=-i,$ 
so  we can adapt the previous approach to see that 
$\|u(t)\|_{L_{x}^2}\to 0$ as $t\to +\infty$. 
However, as we will see below, this is an exceptional case. It turns out 
that highly non-trivial behavior can be observed in \eqref{eq}--\eqref{data} 
for generic small  initial data.

\subsection{The initial value problem for \eqref{eq} with generic small data}

We start the discussion with the following 
basic result.

\begin{thm}[\cite{LNSS1}]\label{thm_initial}
Suppose that $\varphi^0=(\varphi_1^0,\varphi_2^0)\in H^{2}\cap H^{1,1}$ 
and that 
$\|\varphi^0\|_{H^2\cap H^{1,1}}$ is suitably small. 
Let $u=(u_1,u_2)\in C([0,\infty); H^2\cap H^{1,1})$ 
be the solution to \eqref{eq}--\eqref{data}. Then 
there exists $\varphi^+=(\varphi_1^+,\varphi_2^+)\in L^2$ with 
$\hat{\varphi}^+=(\hat{\varphi}_1^+,\hat{\varphi}_2^+)\in L^{\infty}$ 
such that 
\begin{equation}
 \lim_{t\to +\infty} \|u_j(t)-\op{U}(t)\varphi_j^+\|_{L_{x}^2}=0,\quad j=1,2.
\label{scattering_ini}
\end{equation}
Moreover we have 
\begin{equation}
\hat{\varphi}_1^+(\xi)\cdot \hat{\varphi}_2^+(\xi)=0,\quad \xi \in \R.
\label{comp_rel}
\end{equation}
\end{thm}

The global existence part of this asssertion is just a special case of 
more general result (see Theorem~\ref{thm_sdge} below). 
On the other hand, we emphasize that \eqref{comp_rel} should be 
regarded as a consequence of non-trivial long-range nonlinear interactions 
because such a phenomenon does not occur in the usual short-range situation. 
Note also that the system (\ref{eq}) possesses two conservation laws 
\begin{equation*}
 \frac{d}{dt}\Bigl( \|u_1(t)\|_{L_{x}^2}^2+ \|u_2(t)\|_{L_{x}^2}^2 \Bigr)
=
-4\int_{\R} |u_1(t,x)|^2 |u_2(t,x)|^2\, dx
\end{equation*}
and
\begin{equation}
 \frac{d}{dt}\Bigl( \|u_1(t)\|_{L_{x}^2}^2 - \|u_2(t)\|_{L_{x}^2}^2\Bigr)=0.
\label{CL2}
\end{equation}
However, these are not enough to assert that the solution $u=(u_1,u_2)$ is 
asymptotically free in the sense of \eqref{scattering_ini}. 
It is worthwhile to mention that \eqref{CL2} tells us that at least one 
component $u_1(t)$ or $u_2(t)$ does not decay as $t\to +\infty$ in 
$L^2(\R_{x})$ if $\|\varphi_1^0\|_{L^2}\ne \|\varphi_2^0\|_{L^2}$. 
In particular, it is far from obvious whether or not both $u_1(t)$ and 
$u_2(t)$ can behave like non-trivial free solutions as $t\to +\infty$. 
That is why we are interested in (non-)triviality of each component 
of the scattering state.

\subsection{Criteria for (non-)triviality of the scattering state}

To investigate the relation \eqref{comp_rel} in more detail, 
let us point out that we also have the following proposition.

\begin{prop}[\cite{LNSS1}]\label{prop_key_lnss1}
We put $\dis{\varphi_j^+=\lim_{t\to+\infty}\op{U}(-t)u_j(t)}$
in $L^2$, $j=1,2$, for the global solution $u=(u_1,u_2)$ to
\eqref{eq}--\eqref{data},
whose existence is guaranteed by {\upshape{Theorem~ \ref{thm_initial}}}.
There exists a function $m:\R\to\R$ such that the following holds
for each $\xi \in \R$
\begin{itemize}
\item
$m(\xi)>0$ implies
$\hat{\varphi}_1^+(\xi)\ne 0$ and $\hat{\varphi}_2^+(\xi) = 0$\upshape{;}
\item
$m(\xi)<0$ implies
$\hat{\varphi}_1^+(\xi)= 0$ and $\hat{\varphi}_2^+(\xi)\ne 0$\upshape{;}
\item
$m(\xi)=0$ implies
$\hat{\varphi}_1^+(\xi)=\hat{\varphi}_2^+(\xi)=0$.
\end{itemize}
In fact, $m(\xi)$ has the following expression:
\begin{align*}
m(\xi)
=\left|\alpha_1(2,\xi)\right|^2-\left|\alpha_2(2,\xi)\right|^2
+\int_2^{\infty}\rho(\tau,\xi)d\tau,
\end{align*}
where
\begin{align}
\alpha_j(t,\xi)=\op{F} \bigl[\op{U}(-t)u_j(t,\cdot)\bigr](\xi),
\label{def_alpha_j}
\end{align}
\[
 \rho(t,\xi)
 =
 2\realpart\Bigl[ 
  \cc{\alpha_1(t,\xi)}R_1(t,\xi)- \cc{\alpha_2(t,\xi)}R_2(t,\xi) 
  \Bigr],
\]
\begin{align}\label{def_R}
R_1
=
\frac{1}{t} |\alpha_2|^2 \alpha_1
-\op{F}\op{U}(-t)\bigl[ |u_2|^2 u_1\bigr],
\quad
R_2
=
\frac{1}{t} |\alpha_1|^2 \alpha_2
-\op{F}\op{U}(-t)\bigl[ |u_1|^2 u_2\bigr].
\end{align}
\end{prop}

Note that \eqref{comp_rel} follows from Propisition~\ref{prop_key_lnss1} 
immediately. In other words, Propisition~\ref{prop_key_lnss1} is more precise 
than the relation \eqref{comp_rel}, and the function $m(\xi)$ plays an 
important role in it. This indicates that better understanding of $m(\xi)$ 
will bring us more precise information on the scattering state 
$\varphi^+$.
To address this point,  we put a small parameter $\eps$ in front of
the initial data to distinguish information on the amplitude from the others,
that is, we replace the initial condition (\ref{data}) by
\begin{eqnarray}
 u_j(0,x) =\eps \psi_j(x),\quad j=1,2,
\label{data_eps}
\end{eqnarray}
where $\psi_j\in H^{2}\cap H^{1,1}$ is independent of $\eps$. 
Then we have the following.
\begin{thm}[\cite{LNSS2}]\label{thm_m}
Let $m$ be the function given in 
{\upshape{Proposition~\ref{prop_key_lnss1}}} with the initial condition 
\eqref{data} replaced by \eqref{data_eps}. We have
\begin{eqnarray*}
m(\xi)
=
\eps^2 \bigl(|\hat{\psi}_1(\xi)|^2-|\hat{\psi}_2(\xi)|^2\bigr)
+
O(\eps^4)
\end{eqnarray*}
as $\eps \to +0$ uniformly in $\xi \in \R$.
\end{thm}

As a consequence of Theorem~\ref{thm_m}, we have the following criteria
for (non-)triviality of the scattering state 
$\varphi^+=(\varphi_1^+,\varphi_2^+)$ for the initial value problem 
\eqref{eq}--\eqref{data_eps}.

\begin{cor}[\cite{LNSS2}]\label{cor_criterion1}
Assume that there exist points $\xi^* \in \R$ and $\xi_* \in \R$ such that
\begin{eqnarray}
 |\hat{\psi}_1(\xi^*)| > |\hat{\psi}_2(\xi^*)|
\label{m_plus}
\end{eqnarray}
and
\begin{eqnarray}
 |\hat{\psi}_1(\xi_*)| < |\hat{\psi}_2(\xi_*)|,
\label{m_minus}
\end{eqnarray}
respectively.
Then, for sufficiently small $\eps$, we have
$\|\varphi_1^+\|_{L^2}>0$ and $\|\varphi_2^+\|_{L^2}>0$.
\end{cor}

\begin{cor}[\cite{LNSS2}]\label{cor_criterion2}
Assume that
\begin{eqnarray}
 |\hat{\psi}_1(\xi)| > |\hat{\psi}_2(\xi)|
\label{m_plus_everywhere}
\end{eqnarray}
for all $\xi \in \R$. Then, for sufficiently small $\eps$,
$\varphi_2^+$ vanishes almost everywhere on $\R$,
while $\|\varphi_1^+\|_{L^2}>0$.
\end{cor}

It follows from  (\ref{scattering_ini}) and Corollary~\ref{cor_criterion1}
that both $u_1(t)$ and $u_2(t)$ behave like non-trivial free solutions as
$t\to +\infty$. In particular, we see that
$L^2$ decay does not occur for $u_1(t)$ and $u_2(t)$
under (\ref{m_plus}) and (\ref{m_minus}).
To the contrary,  Corollary~\ref{cor_criterion2} tells us
that only the second component $u_2(t)$ is dissipated as $t\to +\infty$ 
in the sense of $L^2$ under (\ref{m_plus_everywhere}).
We emphasize again that such phenomena do not occur in the usual short-range
settings. In this sense, the dynamics for the system (\ref{eq}) is much more
delicate than that for the single Schr\"odinger equation \eqref{nls_1} 
with a dissipative cubic nonlinear term.

At the end of this subsection, let us mention the sketch of the 
proof of Theorem~\ref{thm_m} briefly. 
The key is to focus on the function $\alpha_j$ given by \eqref{def_alpha_j}. 
By the reduction similar to that in the previous section, 
we see that the leading part of $u_j$ as $t\to +\infty$ can be given by 
\[
 \frac{1}{\sqrt{it}}\alpha_j\left(t,\frac{x}{t}\right)e^{i\frac{x^2}{2t}}
\] 
and that the evolution of $\alpha=(\alpha_1,\alpha_2)$ is governed by 
the system
\[
\pa_t \alpha_1 =-\frac{|\alpha_2|^2}{t}\alpha_1 +R_1,\qquad 
\pa_t \alpha_2 =-\frac{|\alpha_1|^2}{t}\alpha_2 +R_2,
\]
where $R_1$ and $R_2$ are given by \eqref{def_R}. If $R_1$ and $R_2$ 
are shown to be harmless, we have
\begin{align}
 \sup_{\xi\in \R}
 \Bigl|
 m(\xi)
 -
 \bigl(|\alpha_1(2,\xi)|^2-|\alpha_2(2,\xi)|^2\bigr)
 \Bigr|
 \le
 C\eps^4.
\label{est_key_1}
\end{align}
Moreover we can show that 
\begin{align}
\alpha_j (2,\xi) 
&=
\hat{u}_j(0,\xi) 
-i\int_0^2 \op{F}\bigl[\op{U}(-t)\op{L}u_j(t,\cdot)\bigr](\xi)\, dt
\notag\\
&=
\eps \hat{\psi}_j(\xi)+ O(\eps^{3})
\label{est_key_2}
\end{align}
as $\eps \to +0$, uniformly in $\xi \in \R$, $j=1$, $2$, provided that 
we have a good control of $u$.  
By \eqref{est_key_1} and \eqref{est_key_2}, we reach the conclusion. 
For the technical details, see \cite{LNSS1} and \cite{LNSS2} 
(see also \cite{NS} and \cite{N} for closely related works on 
the wave equation case.)

\subsection{The final state problem for \eqref{eq}}

To see the role of the relation \eqref{comp_rel} from a different angle, 
let us consider the final state problem for \eqref{eq}, that is, finding a 
solution $u=(u_1,u_2)$ to \eqref{eq} which satisfies
\begin{align}
 \lim_{t\to +\infty} \|u_j(t)-\op{U}(t)\psi_j^+\|_{L_{x}^2}=0,
\qquad j=1,2
\label{scattering}
\end{align}
for a prescribed final state $\psi^+=(\psi_1^+,\psi_2^+)$. Roughly speaking, 
the propositions below imply that 
\eqref{scattering} holds if and only if 
\begin{align}\label{decoupling_fin}
\hat{\psi}_1^+(\xi)\cdot \hat{\psi}_2^+(\xi)=0,\quad \xi \in \R.
\end{align}
Remember that \eqref{scattering} should hold in the short-range case 
regardless of whether \eqref{decoupling_fin} is true or not. In other words, 
our problem must be distinguished from the usual short-range situation. 

The precise statements are as follows.

\begin{prop}[\cite{LNSS1}]\label{prp_final_a}
Let $T_0\geq 1$ be given, and let $u$ be a solution to \eqref{eq} for 
$t\ge T_0$ satisfying 
\begin{align*}
\sup_{t\ge T_0}
\Bigl(
 t^{-\gamma} \|\op{U}(-t)u(t)\|_{H_{x}^{1,1} }
 + 
 \|\op{F}\op{U}(-t)u(t)\|_{L_{x}^{\infty}}
\Bigr)
<\infty
\end{align*}
with some $\gamma\in(0, 1/12)$.
If there exists $\psi^+\in L^2$ with 
$\hat{\psi}^+\in L^{\infty}$ such that \eqref{scattering} holds, 
then we must have \eqref{decoupling_fin}.
\end{prop}

\begin{prop}[\cite{LNSS1}]\label{prp_final_b}
Suppose that  $\psi^+ $ satisfies $\hat{\psi}^+\in H^{0,s}\cap L^{\infty}$ 
with some $s>1$, and that $\delta=\|\hat{\psi}^+\|_{L^{\infty}}$ is suitably small. 
If \eqref{decoupling_fin} holds, then there exist $T\ge 1$ and 
a unique solution $u$ to \eqref{eq} for $t\ge T$ 
satisfying $\op{U}(-t)u\in C([T,\infty);H^{0,1})$
and  \eqref{scattering}. 
\end{prop}

The proof of Proposition~\ref{prp_final_a} is based on a contradiction 
argument. Proposition~\ref{prp_final_b} can be shown by rewriting the system 
\eqref{eq} in the form of integral equations and applying a suitable fixed 
point argument. See \cite{LNSS1} and \cite{Li} for the details of the proof.

\appendix\section{Appendix: General nonlinear Schr\"odinger systems of 
derivative type}
\label{app}

For the convenience of the readers, we collect the results obtained in 
\cite{LiS1} without proof. 
We consider general $n$-component nonlinear Schr\"odinger systems 
in the form 
\begin{align}
\left\{\begin{array}{cl}
\mathcal{L}_{m_j} u_j=N_j(u,\pa_x u), & t>0,\ x\in \R,\ j=1,\ldots, n,\\
u_j(0,x)=\varphi_j(x), & x \in \R,\ j=1,\ldots, n,
\end{array}\right.
\label{nls_system}
\end{align}
where $\mathcal{L}_{m_{j}}=i\pa_t +\frac{1}{2m_{j}}\pa_x^2$, 
$m_j \in \R\backslash \{0\}$, and 
$u=(u_j(t,x))_{1\le j \le n}$ is a $\C^n$-valued unknown function. 
The nonlinear term $N=(N_j)_{1\le j\le n}$ is assumed to be 
a cubic homogeneous polynomial in 
$(u,\pa_x u, \overline{u}, \overline{\pa_x u})$. 
We set $I_n=\{1,\ldots, n\}$ and 
$\sh{I}_n =\{1,\ldots, n, n+1,\ldots, 2n \}$. 
For $z =(z_j)_{j\in I_n}\in \C^n$, we write 
\[
 \sh{z} =(\sh{z}_k)_{k \in \sh{I}_n}
 :=(z_1,\ldots,z_n, \overline{z_1},\ldots, \overline{z_n})
\in \C^{2n}.
\]
Then general cubic nonlinear term $N=(N_j)_{j \in I_n}$ can be written as 
\begin{align*}
 N_j(u,\pa_x u)
 =\sum_{l_1, l_2,l_3=0}^{1}\sum_{k_1, k_2, k_3 \in \sh{I}_n}
 C_{j, k_1, k_2, k_3}^{l_1,l_2,l_3} (\pa_x^{l_1} \sh{u}_{k_1}) 
 (\pa_x^{l_2} \sh{u}_{k_2}) (\pa_x^{l_3} \sh{u}_{k_3})
\end{align*}
with suitable $C_{j, k_1, k_2, k_3}^{l_1,l_2,l_3} \in \C$. 
With this expression of $N$, we define 
$p=(p_j(\xi;Y))_{j \in I_n}:\R \times \C^n \to \C^n$ by  
\begin{align*}
 p_j(\xi; Y)
 :=
 \sum_{l_1, l_2,l_3=0}^{1}\sum_{k_1, k_2, k_3 \in \sh{I}_n}
 C_{j, k_1, k_2, k_3}^{l_1,l_2,l_3} (i\tm_{k_1} \xi)^{l_1}
 (i\tm_{k_2} \xi)^{l_2} (i\tm_{k_3} \xi)^{l_3} 
 \sh{Y}_{k_1} \sh{Y}_{k_2} \sh{Y}_{k_3} 
\end{align*}
for $\xi \in \R$ and $Y=(Y_j)_{j \in I_n} \in \C^n$, where 
\[
 \tm_k=\left\{\begin{array}{cl}
 m_k& (k=1,\ldots, n),\\[4mm]
 -m_{(k-n)} & (k = n+1,\ldots, 2n).
 \end{array}\right.
\]
We denote by $\jb{\cdot, \cdot}_{\C^n}$ the standard scalar 
product in $\C^n$, i.e., 
\[
 \jb{z,w}_{\C^n}=\sum_{j=1}^{n}  z_j \overline{w_j}
\]
for $z=(z_j)_{j \in I_n}$ and $w=(w_j)_{j \in I_n} \in \C^n$. 

With these notations, let us introduce the following conditions: 
\begin{itemize}
\item[(a)] For all $j \in I_n$ and 
$k_1,k_2, k_3 \in \sh{I}_n$, 
\[
 m_j\ne  \tm_{k_1} + \tm_{k_2} + \tm_{k_3}
 \ \ \mbox{implies}\ \ 
 C_{j,k_1,k_2,k_3}^{l_1,l_2,l_3}=0,\ \ 
 l_1, l_2, l_3 \in \{0,1\}.
\]
\item[(b$_0$)] 
There exists an $n\times n$ positive Hermitian matrix $\mathcal{H}$ 
such that 
\[
 \imagpart \jb{p(\xi;Y), \mathcal{H}Y}_{\C^n}\le 0 
\]
for all $(\xi,Y) \in \R\times \C^n$. 
\item[(b$_1$)] 
There exist an $n\times n$ positive Hermitian matrix $\mathcal{H}$ and 
a positive constant $C_*$ such that 
\[
 \imagpart \jb{p(\xi;Y), \mathcal{H}Y}_{\C^n}\le -C_* |Y|^4 
\]
for all $(\xi,Y) \in \R\times \C^n$. 
\item[(b$_2$)] 
There exist an $n\times n$ positive Hermitian matrix $\mathcal{H}$ and 
a positive constant $C_{**}$ such that 
\[
 \imagpart \jb{p(\xi;Y), \mathcal{H}Y}_{\C^n}\le -C_{**} \jb{\xi}^2 |Y|^4
\]
for all $(\xi,Y) \in \R\times \C^n$. 
\item[(b$_3$)] 
$p(\xi;Y)=0$ for all $(\xi,Y) \in \R\times \C^n$. \\
\end{itemize}

We have the following.
\begin{thm}[\cite{LiS1}]  \label{thm_sdge}
Assume the conditions {\upshape{(a)}} and {\upshape{(b$_0$)}} 
are satisfied. Let 
$\varphi=(\varphi_j)_{j \in I_ n} \in H^3\cap H^{2,1}$, 
and assume $\eps:=\|\varphi\|_{H^3}+\|\varphi\|_{H^{2,1}}$ is sufficiently 
small. Then \eqref{nls_system} admits a unique global solution 
$u=(u_j)_{j \in I_n} \in C([0,\infty); H^3\cap H^{2,1})$. 
Moreover we have
\begin{align*}
 \|u(t)\|_{L_{x}^{\infty}} \le \frac{C\eps}{\sqrt{1+t}}, \qquad 
 \|u(t)\|_{L_{x}^2} \le C\eps
\end{align*}
for $t\ge 0$, where $C$ is a positive constant not depending on $\eps$.
\end{thm}

\begin{thm}[\cite{LiS1}, \cite{LNSS3}]  \label{thm_decay1}
Assume the conditions {\upshape{(a)}} and {\upshape{(b$_1$)}} 
are satisfied. Let $u$ 
be the global solution to \eqref{nls_system}, whose existence is guaranteed 
by {\upshape{Theorem \ref{thm_sdge}}}. Then we have 
\[
 \|u(t)\|_{L_{x}^{\infty}} 
\le \frac{C\eps }{\sqrt{(1+t)\{1+ \eps^2\log (2+t)\}}}
\]
for $t\ge 0$, where $C$ is a positive constant not depending on $\eps$. 
We also have
\[
 \|u(t)\|_{L_{x}^{2}} \le \frac{C\eps}{(1+\eps^2 \log (2+t))^{3/8-\delta}},
\]  
where $\delta>0$ can be taken arbitrarily small.
\end{thm}

\begin{thm}[\cite{LiS1}]  \label{thm_decay2}
Assume the conditions {\upshape{(a)}} and {\upshape{(b$_2$)}} are satisfied. 
Let $u$ be as above. Then we have 
\[
 \|u(t)\|_{L_{x}^{2}} \le \frac{C\eps}{\sqrt{1+\eps^2 \log (2+t)}}
\] 
for $t\ge 0$, where $C$ is a positive constant not depending on $\eps$.
\end{thm}

\begin{thm}[\cite{LiS1}]  \label{thm_asymp_free}
Assume the conditions {\upshape{(a)}} and {\upshape{(b$_3$)}} are satisfied. 
Let $u$ be as above. 
For each $j\in I_n$, there exists $\varphi_j^+ \in L^2(\R_{x})$ with 
$\hat{\varphi}_j^+ \in L^{\infty}(\R_{\xi})$ such that 
\[
 u_j(t)
 =
 e^{i\frac{t}{2m_j}\pa_x^2} \varphi_j^+ + O(t^{-1/4+\delta})
 \quad \mbox{in}\ L^2(\R_x)
\]
and 
\[
 u_j(t,x)=\sqrt{\frac{m_j}{i t}}\, \hat{\varphi}^+_j
 \left( \frac{m_j x}{t}\right) e^{i\frac{m_j x^2}{2t}} 
 + O(t^{-3/4+ \delta})
 \quad \mbox{in}\ L^{\infty}(\R_x) 
\]
as $t \to +\infty$, where $\delta>0$ can be taken arbitrarily small.
\end{thm}

In the system case, many interesting problems are left unsolved. 
For more recent progress or related issues, we refer the readers to 
\cite{IKS}, \cite{KatSak}, \cite{KLS}, \cite{Kim}, \cite{KMSU}, \cite{LiS2}, 
\cite{MSU}, \cite{SakSuna}, and so on.

\medskip
\subsection*{Acknowledgments}
The work of Y.~N. is supported by JSPS Grant-in-Aid for JSPS Fellows 
(22J00983). The work of C.~L. is supported by Education
Department of Jilin Province of China (JJKH20220527KJ).
The work of H.~S. is supported by JSPS KAKENHI (21K03314). 
This work is partly supported by 
Osaka Central Advanced Mathematical Institute, Osaka Metropolitan University 
(MEXT Joint Usage/Research Center on Mathematics and Theoretical 
Physics JPMXP0619217849). 


\end{document}